\documentclass[12pt]{article}
\usepackage{amsfonts,amssymb,theorem}

\IfFileExists{mathrsfs.sty}{\usepackage{mathrsfs}\let\mathcal\mathscr}{}

\newcommand{\mystretch}{1.1}

\newcommand{\R}{\mathbb{R}}
\newcommand{\Z}{\mathbb{Z}}
\newcommand{\E}{\mathcal{E}}
\newcommand{\SC}{\mathcal{S}}
\newcommand{\PC}{\mathcal{P}}

\renewcommand{\sp}{\mathfrak{sp}}

\newcommand{\CM}{C^\infty(M)}
\newcommand{\Bbt}{[\![t]\!]}
\newcommand{\CMt}{C^\infty(M)\Bbt}
\def\cyclic{\mathop{\kern0.9ex{{+}
\kern-2.2ex\raise-.28ex\hbox{\Large\hbox{$\circlearrowright$}}}}\limits}
\newcommand{\multiindex}[3]{
\renewcommand{\arraystretch}{0.5}
\begin{array}[t]{c}
#1 \\
{\scriptstyle #2 }\\
{\scriptstyle #3 }
\end{array}
\renewcommand{\arraystretch}{1}
}
\newcommand{\parnot}{:=}
\newcommand{\invariant}{$T^{2n}$-invariant}
\newcommand{\GG}[1]{{\bf #1}}
\newcommand{\n}{\nabla}
\newcommand{\del}{\partial}
\newcommand{\bas}{{}_}
\newcommand{\haut}{{}^}
\makeatletter
\def\operatorname#1{\mathop{\operator@font #1}\nolimits}%
\makeatother

\newcommand{\ad}{\operatorname{ad}}

\newcommand{\id}{\operatorname{id}}

\newcommand{\Tr}{\operatorname{Trace}}

\newtheorem{thm}{Theorem}[section]
\newtheorem{introthm}{Theorem}
\newtheorem{introprop}{Proposition}
\newtheorem{prop}[thm]{Proposition}
\newtheorem{cor}[thm]{Corollary}
\newtheorem{lem}[thm]{Lemma}
\newtheorem{defn}[thm]{Definition}
{\theorembodyfont{\normalfont\rmfamily}

}
\makeatletter
\newenvironment{pf}[1][]%
  {\def\proof@temp{#1}\par\noindent
  \textsc{Proof}\ifx\proof@temp\@empty\else\ (#1)\fi\hspace{1em}}
  {\rule{12pt}{0pt}\hfill\rule{6pt}{6pt}\par\vspace{.4\baselineskip}}
\makeatother

\makeatletter
\def\footnotenomark#1{
  \protected@xdef\@thanks{\@thanks
   \protect\footnotetext{#1}}%
}
\makeatother
\newcommand{\note}[1]{${}^{\hbox{\small(#1)}}$}
\makeatletter
\def\Day{\two@digits\day}\def\Month{\two@digits\month}
\makeatother

\title{Moduli space of symplectic connections of Ricci type\\ on $T^{2n}$ ;
 a formal approach}

\author{M.~Cahen\note{i},\ S.~Gutt\note{i,ii}\footnotenomark{Research of
the first three authors supported by an ARC of the communaut\'e fran\c caise
de Belgique.}\footnotenomark{\textbf{Mathematics Subject
Classification (1991):}\ \ 53C05, 58C35, 53C57.},\\[5pt]
J.~Horowitz\note{i} and J.~Rawnsley\note{iii}
\footnotenomark{(i) Universit\'e Libre de Bruxelles,
Campus Plaine, CP 218,
bvd du triomphe,
1050 Brussels,
Belgium}
\footnotenomark{(ii) Universit\'e de Metz,
Ile du Saulcy,
57045 Metz Cedex 01,
France}
\footnotenomark{(iii) Mathematics Institute,
University of Warwick,
Coventry CV4 7AL,
United Kingdom}
\footnotenomark{Email: mcahen@ulb.ac.be, sgutt@ulb.ac.be
  and j.rawnsley@warwick.ac.uk}
}

\date{~\\[20pt]January 2002}

\renewcommand{\baselinestretch}{\mystretch}

\begin{document}

\setcounter{page}{0}
\renewcommand{\baselinestretch}{1} 
\thispagestyle{empty}
\maketitle

\thispagestyle{empty}

\begin{abstract}
We consider analytic curves $\n^t$ of symplectic connections of Ricci
type on the torus $T^{2n}$ with $\n^0$ the standard connection. We show,
by a recursion argument, that if $\n^t$ is a formal curve of such
connections then there exists a formal curve of symplectomorphisms
$\psi_t$ such that $\psi_t\cdot\n^t$ is a formal curve of flat
\invariant\ symplectic connections and so $\n^t$ is flat for all $t$.
Applying this result to the Taylor series of the analytic curve, it
means that analytic curves of symplectic connections of Ricci type
starting at $\n^0$ are also flat.

The group $G$ of symplectomorphisms of the torus $(T^{2n},\omega)$ acts
on the space $\E$ of symplectic connections which are of Ricci type. As
a preliminary to studying the moduli space $\E/G$ we study the moduli of
formal curves of connections under the action of formal curves of
symplectomorphisms.
\end{abstract}

\newpage

\renewcommand{\baselinestretch}{\mystretch}

\section*{Introduction~~}

\noindent On any symplectic manifold $(M,\omega)$ the space $\SC$ of
symplectic connections is an infinite dimensional affine space whose
corresponding vector space is the space of completely symmetric
$3$-tensors on $M$. To encode some geometry into a symplectic connection
it thus seems reasonable to introduce a selection rule for symplectic
connections. A variational principle associated to a Lagrangian density,
which is an invariant quadratic polynomial in the curvature, has been
considered in \cite{bib:BC}; the symplectic connections satisfying  the
Euler--Lagrange equations are said to be \textit{preferred}. The
symplectomorphism group $G$ of $(M,\omega)$ acts naturally on $\SC$ and
stabilises the subspace $\PC$ of preferred symplectic connections. The
first question we wanted to address is to give a description of the
moduli space $\PC/G$ of preferred connections modulo the action of
symplectomorphisms. Such a description was given in \cite{bib:BC} when
$(M,\omega)$ is a closed surface; but, up to now, very little has been
done in the higher dimensional situation.

We have observed that a linear condition on the curvature (the vanishing
of one of its irreducible components -- the non-Ricci component, $W$)
implies the Euler--Lagrange equations. Furthermore, this condition seems
to imply that many of the properties of the surface situation extend to
the higher-dimensional case. We have called symplectic connections
satisfying this curvature condition \textit{connections of Ricci type}
(all symplectic connections in dimension $2$ are of Ricci type). This
condition is preserved by symplectomorphisms and so we modify our
initial question to the following one: give a description of the space
$\E$ of Ricci type connections and its moduli space $\E/G$.

This paper is devoted to this modified question in the case where $M$ is
a torus $T^{2n}$ and $\omega$ a  $T^{2n}$-invariant symplectic structure.
 Although we do
not answer this question, we are able, in a formal setting made precise
below, to show that the moduli space is infinite dimensional and to give a
partial description of it.

If $\n^t$ is a formal curve of symplectic connections, we shall denote
by $W^t$ the $W$ part of the curvature of $\n^t$. We prove

\begin{introthm}
Let $\n^t$ be a formal curve of symplectic connections on
$(T^{2n},\omega)$ such that $\n^0$ is the standard flat connection
on $T^{2n}$, and such that $W^t=0$. Then
the formal curvature $R^t$ of $\n^t$ vanishes and there
exists a formal curve of symplectomorphisms $\psi_t$ such that
${\widetilde\n}^t \parnot \psi_t.\n^t$ is a formal curve of
flat \invariant\ symplectic connections.
\end{introthm}

\noindent This implies

\begin{introthm}
Let $\n^t$ be an analytic curve of analytic symplectic connections on
$(T^{2n},\omega)$ such that $\n^0$ is the standard flat connection
on $T^{2n}$, and such that $W^t=0$. Then
the  curvature $R^t$ of $\n^t$ vanishes.
\end{introthm}

For the moduli space in the formal setting, we show:
\begin{introprop}
For two curves ${\widetilde{\n^t}}$ and ${\widetilde{{\n'}^t}}$
  of invariant flat connections
of Ricci-type on $(\R^{2n}, \Omega)$ with
${\widetilde{\n^0}}={\widetilde{{\n'}^0}}$ the trivial connection,
there always exists a formal curve of symplectomorphisms
${\widetilde{\psi_t}}$ so that ${\widetilde{\psi_t}}
\cdot {\widetilde{\n^t}}={\widetilde{{\n'}^t}}$.
\end{introprop}

\begin{introthm}
The moduli space of formal curves of Ricci-type symplectic connections
starting with the standard flat connection on $(T^{2n},\omega)$ under
the action of formal curves of symplectomorphisms is described by the
space of formal curves of linear maps $A^t \colon \R^{2n} \to
\sp(2n,\R)\Bbt$ satisfying $A^t(X)A^t(Y)=0$ and $A^t(X)Y=A^t(Y)X$,
modulo the action of $Sp(2n,\Z)$.
\end{introthm}

The plan of the paper is as follows. In $\S$\ref{section:curvature} we
recall some general properties of symplectic connections having
Ricci-type curvature.
In $\S$\ref{section:formal} we introduce the notion of formal curves
of connections  and we show that the properties of
$\S$\ref{section:curvature} are still true for a formal curve of
symplectic connections with Ricci-type curvature.
In $\S$\ref{section:torus}, we analyse the $W^t=0$ condition at order $1$
and order $2$ for $\n^{t}=\n^0+\sum_{k=1}^\infty t^kA^{(k)}$   a formal
curve of Ricci-type symplectic connections on $T^{2n}$ with $\n^0$  the
standard flat connection; in particular, we show that there exists a
function $U^{(1)}$ and a completely symmetric, \invariant\ 3-tensor
$Q^{(1)}$ on $T^{2n}$ such that ${\underline{A}}^{(1)} =
({\n}^0)^3U^{(1)}+ Q^{(1)}$ and we show that
$\n^{'t}=\n^0+t{\overline{Q}}^{(1)}$ (with
$\omega({\overline{Q}}^{(1)}(X)Y,Z)=Q^{(1)}(X,Y,Z)$) defines a curve of
invariant flat symplectic connections on $(T^{2n},\omega)$. This remark
can be formulated in a slightly different way: given $\n^t=\n^0+A^{(t)}$
  a smooth curve of Ricci-type symplectic  connections then, up to a
symplectomorphism, the tangent vector to this family of connections lies
in the finite dimensional space of flat \invariant\ symplectic
connections.
$\S$\ref{section:recurrence} is devoted to a proof of a
recurrence lemma which implies the first theorem.
In $\S$\ref{section:eq} we study the question of when two formal curves of flat
invariant connections on $T^{2n}$ are equivalent by a formal curve of
symplectomorphisms.

\smallskip

\noindent{\bf Thanks: }{ We would like to thank Boguslaw Hajduk and Aleksy
Tralle
who pointed out a mistake in an earlier version of this paper.}

\section{Ricci Type Curvature}\label{section:curvature}

 A symplectic connection $\n$  on a symplectic manifold
$(M,\omega)$ is
a linear connection having no torsion and for which $\omega$ is parallel
($\n\omega=0$).
The curvature endomorphism $R$ of $\n$ is
\[
R(X,Y)Z = \left(\nabla_X\nabla_Y - \nabla_Y\nabla_X -
\nabla_{[X,Y]}\right)Z
\]
for vector fields $X,Y,Z$ on $M$. The symplectic curvature tensor
\[
R(X,Y;Z,T) = \omega(R(X,Y)Z,T)
\]
is antisymmetric in its first two arguments, symmetric in its last two
and satisfies the first Bianchi identity
\[
\cyclic_{X,Y,Z}\ R(X,Y;Z,T) = 0
\]
where $\cyclic$ denotes the sum over the cyclic permutations of the listed
set of elements. The second Bianchi identity takes the form
\[
\cyclic_{X,Y,Z} \ \left(\nabla_XR\right)(Y,Z) = 0.
\]
The Ricci tensor $r$ is the symmetric $2$-tensor
\[
r(X,Y) = \Tr[ Z \mapsto R(X,Z)Y].
\]
If $\dim M=2n\ge 4$, the curvature $R$ of such a connection has
2 irreducible components under the action of the symplectic group
$Sp(2n,\R)$. We denote them by $E$ and $W$:
\[
R = E + W.
\]
The $E$ component encodes the information contained in the Ricci tensor
of $\n$ and is called the Ricci part  of the curvature tensor. It is given by
\begin{eqnarray*}
E(X,Y;Z,T) &=& \frac{-1}{2(n+1)} \biggl[2\omega(X,Y)r(Z,T) +
\omega(X,Z)r(Y,T) + \omega(X,T)r(Y,Z)\nonumber\\
&&\qquad\mbox{}-\omega(Y,Z)r(X,T) - \omega(Y,T)r(X,Z)
\biggr].
\end{eqnarray*}
The curvature is said to be \textit{of Ricci type} if the $W$
component vanishes, i.e. when $R=E$.

\begin{lem}\label{l1}
Let $(M,\omega)$ be a symplectic manifold of dimension $2n\ge4$. If the
curvature of a symplectic connection $\nabla$ on $M$ is of Ricci type
then there is a $1$-form $u$ such that
\[
\left(\nabla_Xr\right)(Y,Z) = \frac{1}{2n+1}\left(\omega(X,Y)u(Z)
+ \omega(X,Z)u(Y)\right).
\]
Conversely, if there is such a $1$-form $u$, the ``Weyl'' part of the
curvature, $W=R-E$ satisfies
\[
\cyclic_{X,Y,Z} \left(\nabla_X W\right)(Y,Z;T,U) = 0.
\]
\end{lem}

\begin{pf}
The property follows from the second Bianchi's identity, see \cite{bib:CGHR}.
\end{pf}

\begin{cor}
A symplectic manifold with a symplectic connection whose curvature is of
Ricci type is locally symmetric if and only if the $1$-form $u$, defined
in the lemma, vanishes.
\end{cor}
Denote by $\rho$ the linear endomorphism such that
\[
r(X,Y) = \omega(X, \rho Y).
\]
The symmetry of $r$ is equivalent to saying that
$\rho$ is in the Lie algebra of the symplectic group $Sp(TM,\omega)$. For
an integer $p>1$, define
\[
\stackrel{(p)}r(X,Y)=\omega(X,\rho^pY).
\]
It is symmetric when $p$ is odd and antisymmetric when $p$ is even.

\begin{lem}\label{l2}
Let $(M,\omega)$ be a symplectic manifold with a symplectic connection
$\nabla$ with Ricci-type curvature. Then, the following identities hold:
\begin{enumerate}
\item There is a function $b$ such that
\[
\nabla u = -\frac{1+2n}{2(1+n)}\stackrel{(2)}{r} + b\omega.
\]
\item The differential of the function $b$ is given by
\[
db = \frac{1}{1+n} i(\overline{u})r
\]
where $\overline u$ is the vector field such that $i(\overline
u)\omega=u$, so that
\item
\[
b+ \frac{2n+1}{4(1+n)} \Tr \rho^2
\]
is a constant when $M$ is connected.
\end{enumerate}
\end{lem}

\begin{pf}
These identities follow from Lemma \ref{l1}, see \cite{bib:CGHR}.
\end{pf}

 Let the torus $T^{2n}$ be endowed with a  $T^{2n}$-invariant symplectic
structure $\omega$.  Let $\n$ be a symplectic connection
on $(T^{2n},\omega)$ which is of Ricci-type.
The group $G$ of symplectomorphisms of $(T^{2n},\omega)$ acts on
the set $\E$ of symplectic connections with $W=0$. We are interested in
the set of orbits of G in $\E$, i.e.\ in $\E/G$.

\medskip

We now consider the symplectic vector space $(\R^{2n},\Omega)$ and view
$\Omega$ as a translation invariant symplectic structure. A symplectic
connection on $\R^{2n}$ will be determined by its values on translation
invariant vector fields. If, in addition, the connection $\n$ is
translation invariant then $B(X)Y \parnot \n_XY$ (for invariant vector
fields $X,Y$) defines a linear map $B \colon \R^{2n} \to \sp(2n,\R)$
which completely determines $\n$. The only condition on $B$ is that
$\Omega(B(X)Y,Z)$ is completely symmetric.

\begin{prop}\label{invariantRiccitype=flat}
Let $\n$ be a translation invariant symplectic connection on
$(\R^{2n},\Omega)$ and let $B(X)Y = \n_XY$ as above. If $\n$ is of Ricci
type and $2n \ge 4$, then $\n$ is flat and $B(X)B(Y)=0$.
\end{prop}

\begin{pf}
Since $B$ is constant, the curvature endomorphism is given by
\[
R(X,Y) = [B(X),B(Y)]
\]
and so the Ricci tensor is given by
\[
r(X,Y) = \Tr(B(X)B(Y)).
\]
It is easy to see that symplectic curvature tensors $R(X,Y;Z,T)$ are, in
fact, determined by the terms of the form $R(X,Y;X,Y)$ so that the
equation $W=0$ is equivalent to $R(X,Y;X,Y) =
-\frac{2}{n+1}\Omega(X,Y)r(X,Y)$, and in the present case this has the
form
\[
(n+1) \Omega(B(X)X,B(Y)Y) = - 2 \Omega(X,Y) r(X,Y).
\]
Polarising the equation in $X$ we have
\begin{eqnarray*}
(n+1) \Omega(T,B(X)B(Y)Y)
&=&  \Omega(X,Y) r(T,Y) + \Omega(T,Y) r(X,Y)\\
&=&  \Omega(X,Y) \Omega(T,\rho Y) + \Omega(T,Y) \Omega(X,\rho Y),
\end{eqnarray*}
so that $W=0$ is equivalent to
\[
(n+1) B(X)B(Y)Y = \Omega(X,Y) \rho Y + \Omega(X,\rho Y) Y.
\]
Polarising this in $Y$ we have
\[
2(n+1) B(X)B(Y)Z = \Omega(X,Y) \rho Z + \Omega(X,\rho Y)Z + \Omega(X,Z)
\rho Y +
\Omega(X,\rho Z) Y \quad(*).
\]
Now choose dual bases $X^i$, $X_i$ for $\R^{2n}$ with $\Omega(X^i,X_j) =
\delta^i_j$ then an easy calculation shows
\[
\rho  = \sum_i B(X^i) B(X_i).
\]
If we multiply $(*)$ by $B(X^i)$, set $X=X_i$ and sum we get
\[
(n+1) \rho  B(Y)Z = -  B(Y) \rho  Z -  B(Z) \rho  Y.
\]
Alternatively we may substitute $B(X_i)Z$ for $Z$ in (*), set $Y=X^i$
and sum to give
\[
(n+1) B(X) \rho  Z = -  \rho  B(Z) X + B(Z) \rho  X.
\]
Adding the two equations after setting $X=Y$ we see that
\[
\rho  B(X) = - B(X) \rho
\]
and hence that
\[
(n-1) \rho  B(X) = 0.
\]
Thus if $2n \ge 4$
\[
\rho B(X) =  B(X) \rho = 0 \quad \Rightarrow \rho^2 = 0.
\]
Substituting $\rho Z$ for $Z$ in (*) we have
\[
0 = r(X,Y) \rho Z + r(X,Z)\rho Y
\]
and setting $Z=Y$, applying $\Omega(X, \,.\,)$ we get finally
\[
0 = r(X,Y)^2.
\]
Thus the Ricci tensor vanishes, and hence $\n$ is flat.

Putting $\rho = 0$ in (*) yields $B(X)B(Y) = 0$.
\end{pf}

\section{Formal curves}\label{section:formal}

\begin{defn}
A \GG{formal curve of symplectic connections} on a symplectic manifold
$(M,\omega)$ is a formal power series
\[
\n^t = \n + \sum_{k=1}^\infty t^k A^{(k)}
\]
where $\n$ is a symplectic connection on $M$, and the $A^{(k)}$ are
$(2,1)$ tensors such that
\begin{equation}\label{Abar}
{\underline{A}^{(k)}}(X,Y,Z)\colon =\omega(A^{(k)}(X)Y,Z)
\end{equation}
is totally symmetric.
\end{defn}

\begin{defn}
A \GG{formal curve of symplectomorphisms} is a
homomorphism of Poisson algebras
\[
\psi_t \colon \CM \longrightarrow \CMt,
\qquad \psi_t = \psi^{(0)}+
\sum_{k=1}^\infty t^k\psi^{(k)}
\]
such that $\psi^{(0)} \colon \CM \longrightarrow \CM$ is an
isomorphism.
\end{defn}

The leading term $\psi^{(0)}$ of a formal curve of symplectomorphisms is
given by composition with a symplectomorphism $\psi^{(0)}(f) = f \circ
\sigma = \sigma^*(f)$ so that we may take such a term out as a common
factor and write $\psi_t = \sigma^* \circ \phi_t$ and $\phi_t = \id +
\sum_{k\ge1} t^k \phi^{(k)}$.

If $\phi_t = \id + \sum_{k\ge1} t^k \phi^{(k)}$ is a formal curve of
symplectomorphisms beginning with the identity then the first order term
$X^{(1)}=\phi^{(1)}$ is a symplectic vector field. Moreover, for any
symplectic vector field, $\exp tX = \id + \sum_{k \ge 1} t^k/k! X^k$ is
a formal curve of symplectomorphisms. A straightforward recursion
argument then shows that any formal curve of symplectomorphisms
beginning with the identity can be written in the form $\phi_t = \exp
X_t$ where $X_t = \sum_{k \ge 1} t^k X^{(k)}$ is a formal curve of
vector fields.

\begin{defn}
A \GG{formal $1$-parameter group of symplectomorphisms} is a formal
curve of symplectomorphisms $\psi_t$ such that $\psi_{at} \circ
\psi_{bt} = \psi_{(a+b)t}$ for all $a,b \in \R$.
\end{defn}

In order for this definition to make sense we first have to extend
$\psi_t$ by linearity over $\R\Bbt$ to a morphism of $\R\Bbt$-algebras.
The definition then implies that $\psi^{(0)}$ is the identity and that
$\psi^{(1)}(f) = X(f)$ for some symplectic vector field which we call
\GG{the infinitesimal generator} of $\psi_t$. It is easy to see that
every formal $1$-parameter group of symplectomorphisms has the form
$\psi_t = \exp tX$. Moreover, a recursion shows that, if $X_t$
is a formal curve of symplectic vector fields, we can find a
second sequence of symplectic vector fields $Y^{(k)}$ such that
\[
\exp X_t = \exp tY^{(1)} \circ \exp t^2 Y^{(2)} \circ \cdots \circ \exp
t^k Y^{(k)} \circ \cdots
\]
and so any formal curve of symplectomorphisms $\psi_t$ can be factorised
in two ways
\[
\psi_t = \sigma^* \circ \exp X_t = \sigma^* \circ \phi_t^{(1)} \circ
\phi_{t^2}^{(2)} \circ \cdots \circ \phi_{t^k}^{(k)} \circ \cdots
\]
where the the $\phi_{t}^{(k)}$ are formal 1-parameter groups of
symplectomorphisms.

\medskip

Remark that a formal curve of symplectomorphisms $\psi_t$ acts
on a formal curve of vector fields $X_t$ viewed as a $\R\Bbt$-linear
derivation of $\CMt$ by
\[
(\psi_t\cdot X_t)f=\psi_t(X_t({\psi_t}^{-1} f)),
\]
and acts on a formal curve of symplectic connections $\n^t$ by
\begin{equation}\label{actconn}
(\psi_t\cdot \n^t)_X Y= \psi_t\cdot \left(\n^t_{{\psi_t}^{-1}
\cdot X}{\psi_t}^{-1}\cdot Y\right).
\end{equation}

Let $\n^t$ be a formal curve of symplectic connections on a symplectic
manifold $(M,\omega)$ of dimension $2n$,
\[
\n^{t}=\n+\sum_{k=1}^\infty t^kA^{(k)}.
\]
  We denote as in (\ref{Abar}) by
${\underline{A}}^{(k)}$
the corresponding symmetric $3$-tensors.
The formal curvature endomorphism $R^t$ of $\n^t$ is $
R^t(X,Y) = \nabla^t_X \circ\nabla^t_Y - \nabla^t_Y\circ\nabla^t_X -
\nabla^t_{[X,Y]}$ so that
\[
R^{t}=R^{\n}+\sum_{k=1}^\infty t^kR^{(k)}
\]
with
\begin{equation}\label{defR}
R^{(k)}(X,Y)=(\n_X A^{(k)})(Y)-(\n_Y A^{(k)})(X)+
{\multiindex{\sum}{p+q=k}{p,q\ge 1}}[A^{(p)}(X),A^{(q)}(Y)].
\end{equation}
The symplectic curvature tensor $R^t(X,Y;Z,T) = \omega(R^t(X,Y)Z,T)$
is antisymmetric in its first two arguments, symmetric in its last two,
satisfies the first Bianchi identity $\cyclic_{X,Y,Z}\ R^t(X,Y;Z,T) = 0$
and the second Bianchi identity  $\cyclic_{X,Y,Z}
  \ \left(\nabla^t_XR^t\right)(Y,Z) = 0.$

  The formal Ricci tensor is $r^t(X,Y) = \Tr[ Z \mapsto R^t(X,Z)Y]$, so that
\[
r^{t}=r^{\n}+\sum_{k=1}^\infty t^kr^{(k)}
\]
where the $r^{(k)}$ are the symmetric tensors
\begin{equation}\label{defr}
r^{(k)}(X,Y)=\Tr[ Z \mapsto(\nabla_ZA^{(k)})(X)Y]+
{\multiindex{\sum}{p+q=k}{p,q\ge 1}}\Tr A^{(p)}(X)A^{(q)}(Y).
\end{equation}
The Ricci part $E^t$ of the formal curvature tensor is given by
\begin{eqnarray}\label{wt=0}
E^t(X,Y;Z,T) &=& \frac{-1}{2(n+1)} \biggl[2\omega(X,Y)r^t(Z,T) +
\omega(X,Z)r^t(Y,T) + \omega(X,T)r^t(Y,Z)\nonumber\\
&&\qquad\mbox{}-\omega(Y,Z)r^t(X,T) - \omega(Y,T)r^t(X,Z)
\biggr].
\end{eqnarray}
The formal curvature is said to be \textit{of Ricci type} when $R^t=E^t$.
\begin{lem}\label{l1formal}
Let $(M,\omega)$ be a symplectic manifold of dimension $2n\ge4$. If the
formal curvature of a formal curve of symplectic connections $\nabla^t$
on $M$ is of Ricci type then there exists a formal curve of $1$-forms
\[
u^{t} = \sum_{k=0}^\infty t^k u^{(k)}
\]
such that
\begin{equation} \label{eqderRicci}
(\n^t_X r^t)(Y,Z) =
\frac1{2n+1}(\omega(X,Y)u^t(Z)+\omega(X,Z)u^t(Y))
\end{equation}
and there exists a formal curve of functions
\[
  b^t = \sum_{k=0}^\infty t^k b^{(k)}
\]
such that
\begin{equation} \label{eqderu}
  \n^t u^t = -\frac{1+2n}{2(1+n)}\stackrel{(2)}{r^t} + b^t\omega.
\end{equation}
   with $\omega(X,(\rho^t)Y)=r^t(X,Y)=$ and
$\stackrel{(2)}{r^t}(X,Y)=\omega(X,(\rho^t)^2Y)$.
Also
\begin{equation} \label{eqderb}
  db^t = \frac1{1+n}i(\overline u^t)r^t.
\end{equation}
\end{lem}

\begin{lem}\label{flatformal}
Let $\n^t$ be a formal curve of translation invariant symplectic
connections on $(\R^{2n},\Omega)$ and let $B^t(X)Y  \parnot \n^t_XY$
(for invariant vector fields $X,Y$). If $\n^t$ is of Ricci type and $2n
\ge 4$, then $\n^t$ is flat and $B^t(X)B^t(Y)=0$.
\end{lem}

\begin{pf}
We can copy in the formal series setting the proof
of Lemma \ref{invariantRiccitype=flat}.
Write $B^t=\sum_{k=0}^{\infty} t^kB^{(k)}$ where the $B^{(k)}$ are
constant maps from $\R^{2n}$ to $sp(\R^{2n},\Omega)$.
The formal curvature endomorphism is given by
\[
R^t(X,Y) = [B^t(X),B^t(Y)] \quad \mbox{i.e.}\quad
R^{(k)}(X,Y) = {\multiindex{\sum}{p+q=k}{p,q\ge 0}}[B^p(X),B^q(Y)]
\]
and the formal Ricci tensor  by
\[
r^t(X,Y) = \Tr(B^t(X)B^t(Y)) \quad \mbox{i.e.}\quad
r^{(k)}(X,Y) = {\multiindex{\sum}{p+q=k}{p,q\ge 0}}\Tr B^p(X) B^q(Y) .
\]
The equation $W^t=0$ is again equivalent to $
2(n+1) B^t(X)B^t(Y)Z = \Omega(X,Y) \rho^tZ + \Omega(X,\rho^tY)Z +
\Omega(X,Z) \rho^tY +
\Omega(X,\rho^tZ) Y $, i.e.
\begin{eqnarray}\label{*k}
{\multiindex{\sum}{p+q=k}{p,q\ge 0}}2(n+1) B^{(p)}(X)B^{(q)}(Y)Z
&=&\Omega(X,Y) \rho^{(k)}Z
+ \Omega(X,\rho^{(k)}Y)Z \nonumber  \\
& &+ \Omega(X,Z) \rho^{(k)}Y +
\Omega(X,\rho^{(k)}Z) Y.
\end{eqnarray}
Choosing dual bases $X^i$, $X_i$ for $\R^{2n}$ with $\Omega(X^i,X_j) =
\delta^i_j$ then $\rho^t  = \sum_i B^t(X^i) B^t(X_i)$, i.e.
$\rho^{(k)} = \sum_{p+q=k}\sum_i B^{(p)}(X^i) B^{(q)}(X_i).$
If we multiply (\ref{*k}) by $B^{(k')}(X^i)$, set $X=X_i$ and sum over $i$
and over $k,k'\ge 0$ so that $k+k'=K$ we get
\[
(n+1) {\multiindex{\sum}{q'+q=K}{q,q'\ge 0}}\rho^{(q')} B^{(q)}(Y)Z =
       {\multiindex{\sum}{k'+k=K}{k',k'\ge 0}}
       \left( -  B^{(k')}(Y) \rho ^{(k)} Z -  B^{(k')}(Z) \rho^{(k)}  Y\right).
\]
This can be written in terms of formal series
\[
(n+1) \rho^t B^t(Y)  Z = -   B^t(Y)  \rho^t Z - B^t(Z) \rho^t Y.
\]

Alternatively we may substitute $B^{(s)}(X_i)Z$ for $Z$ in (\ref{*k}), set
$Y=X^i$
and sum to give
\[
(n+1) B^t(X) \rho^t  Z = -  \rho^t  B^t(Z) X + B^t(Z) \rho^t X.
\]
Adding the two equations after setting $X=Y$ as before, we see that
$\rho^t  B^t(X) = - B^t(X) \rho^t $,
so $(n-1) \rho^t  B^t(X) = 0$ and, if $2n \ge 4$,
$\rho^t B^t(X) =  B^t(X) \rho^t = 0$ thus $(\rho^t)^2 = 0$.
This in turn implies $r^t=0$, hence $R^t=0$ and
  $\n$ is flat.
Putting $\rho^t = 0$ in \ref{*k} yields $B^t(X)B^t(Y) = 0$.
\end{pf}

\section{Curves of Ricci Type Connections on the Torus}\label{section:torus}
Consider the torus $T^{2n}$ endowed with a \invariant\ symplectic
structure $\omega$.
Let $\n^0$ be the standard flat, \invariant\
symplectic connection on $(T^{2n},\omega)$.
Let
\[
\n^{t}=\n^0+\sum_{k=1}^\infty t^kA^{(k)}
\]
  be a formal curve of
symplectic connections such that $W(t)=0$. We denote as before (\ref{Abar})
by ${\underline{A}}^{(k)}$
the corresponding symmetric $3$-tensors
(${\underline{A}}^{(k)} (X,Y,Z) = \omega(A^{(k)}(X)Y,Z)$).

We consider, as given by Lemma \ref{l1formal}, the corresponding
formal curve of $1$-forms
$u^{t} = \sum_{k=0}^\infty t^k u^{(k)}$
  and the  formal curve of functions
  $b^t = \sum_{k=0}^\infty t^k b^{(k)}$; clearly
$u^{(0)}=0$
and $b^{(0)} = 0$ since $r^{\n^0}=0$.

\begin{lem}\label{order1}
If $\n^{t}=\n^0+\sum_{k=1}^\infty t^kA^{(k)}$ is a formal curve of
symplectic connections such that $W(t)=0$, then the formal curvature
vanishes at order $1$ in $t$ (i.e. one has $b^{(1)}=0$, $u^{(1)}=0$,
$r^{(1)}=0$, $R^{(1)}=0$). Furthermore, there exists a function
$U^{(1)}$ and a completely symmetric, \invariant\ 3-tensor $Q^{(1)}$ on
$T^{2n}$ such that
\[
{\underline{A}}^{(1)} = ({\n}^0)^3U^{(1)}+ Q^{(1)}.
\]
\end{lem}

\begin{pf}
Denote by $x^a$ ($1 \le a\le 2n$) the standard angle
variables on $T^{2n}$ and by $\del_a$ the corresponding \invariant\
vector fields on $T^{2n}$ (the standard flat connection is defined by
$\n^0_{\del_a}\del_b = 0$).

At order 1, since $b^{(0)} = 0$, $u^{(0)}=0$, $r^0=0$, we have:
\begin{enumerate}
\item $db^{(1)}=0$ by (\ref{eqderb}), so $b^{(1)}$ is a constant;
\item $du^{(1)}=b^{(1)}\omega$ by (\ref{eqderu}); but $\omega$ is not exact
  by compactness of $T^{2n}$
so $b^{(1)}=0$ and $\n^0 u^{(1)}=0$ thus $u^{(1)}(X)$ is a constant for any
\invariant vector field $X$ on $T^{2n}$;
\item the equation (\ref{eqderRicci}) at order $1$ yields
$(\n^0 r^1)$ as a combination
of products of $\omega$ and $u^1$ so that
$\del_a(r^{(1)}(\del_b,\del_c))$ is a constant;
the periodicity of the angles $x^a$ implies then that
$\del_a(r^{(1)}(\del_b,\del_c))=0$
so $u^{(1)}=0$ and $r^{(1)}(\del_b,\del_c)=a^{(1)}_{ab}$ is a constant.
\end{enumerate}

The definition of the (formal) Ricci tensor (\ref{defr}) at order $1$
yields $a^{(1)}_{ab} = -\del_q {A^{(1)}}^q_{ab}$; hence, for each value
of the indices $a,b$, the $2n$-form $a^{(1)}_{ab}\omega^n$ is exact;
this implies
\[
  a^{(1)}_{ab}=0\quad ~\mbox{so}\quad ~ r^{(1)}=0\quad ~
  \mbox{and thus}\quad ~ R^{(1)} = 0.
\]
The definition of the (formal) curvature tensor (\ref{defR}) at order
$1$ gives $R^{(1)}_{abcd}=\del_a {\underline{A}}^{(1)}_{bcd} -\del_b
{\underline{A}}^{(1)}_{acd}$. Hence, for each value of the indices $c,d$
the 1-form ${\underline{A}}^{(1)}_{. cd}$ is closed, so there exist
functions $k_{cd}$ on $T^{2n}$ and constants $Q^{(1)}_{bcd}$ such that:
\[
  {\underline{A}}^{(1)}_{bcd} = \del_b k^{(1)}_{cd}+ Q^{(1)}_{bcd}.
\]
Since $\n^t$ is symplectic, ${\underline{A}}^{(1)}_{bcd}$ is totally
symmetric; the fact that
${\underline{A}}^{(1)}_{bcd}-{\underline{A}}^{(1)}_{cbd} = 0$ implies
\[
  \del_b k^{(1)}_{cd} - \del_c k^{(1)}_{bd}= - Q^{(1)}_{bcd}
+ Q^{(1)}_{cbd}.
\]
When $d$ is fixed, the left-hand side is an exact 2-form. The right-hand
side is \invariant. Since there are no non-zero exact \invariant\ forms,
this implies
\[
  Q^{(1)}_{bcd} = Q^{(1)}_{cbd}, \qquad \qquad
\del_b k^{(1)}_{cd} - \del_c k^{(1)}_{bd} = 0.
\]
Similarly ${\underline{A}}^{(1)}_{bcd}-{\underline{A}}^{(1)}_{bdc} = 0$ gives
\[
  \del_b k^{(1)}_{cd} - \del_b k^{(1)}_{dc} = - Q^{(1)}_{bcd} +
- Q^{(1)}_{bdc}.
\]
In this case, when $c$ and $d$ are fixed, the left-hand side is an exact
1-form, while the right-hand side is \invariant. For the same reason as
above, we deduce that both members vanish:
\[
  Q^{(1)}_{bcd} = Q^{(1)}_{bdc} \qquad \qquad k^{(1)}_{cd}-
k^{(1)}_{dc}=\mbox{constant}.
\]
Hence $Q^{(1)}_{bcd}$ is completely symmetric. Furthermore, for each
fixed index $d$, the 1-form $k^{(1)}_{. d}$ is closed. Hence there exist
functions $S^{(1)}_d$ and constants $T_{cd}$ such that
\[
  k^{(1)}_{cd}=\del_c S^{(1)}_d + T^{(1)}_{cd}.
\]
The fact that $k^{(1)}_{cd}-k^{(1)}_{dc}$ is a constant implies
for the 1-form $S^{(1)}_{.}$  that $dS^{(1)}$ is \invariant,
thus $S^{(1)}$ is closed. Hence there exists a function $U^{(1)}$
and constants $V^{(1)}_d$  such that
\[
  S^{(1)}_d = \del_d U^{(1)}+V^{(1)}_d.
\]
  Substituting, we have:
\[
{\underline{A}}^{(1)}_{bcd} = \del^3_{bcd}U^{(1)}+ Q^{(1)}_{bcd}.
\]
\end{pf}

\begin{lem}\label{order2}
If $\n^{t}=\n^0+\sum_{k=1}^\infty t^kA^{(k)}$
is a formal curve of
symplectic connections such that $W(t)=0$, then
the curvature vanishes at order $2$ in $t$,
(i.e. $b^{(2)}=0$, $u^{(2)}=0$, $r^{(2)}=0$, $R^{(2)}=0$).

Writing ${\underline{A}}^{(1)} = ({\n}^0)^3U^{(1)}+ Q^{(1)}$
as in Lemma \ref{order1},
the formula
$\n^{'t}=\n^0+t{\overline{Q}}^{(1)}$, where
$\omega({\overline{Q}}^{(1)}(X)Y,Z)=Q^{(1)}(X,Y,Z)$,
defines a curve of invariant flat symplectic connections on $(T^{2n},\omega)$.

Furthermore, there exist a function $U^{(2)}$ and a \invariant,
completely symmetric tensor $Q^{(2)}$ such that
\[
{\underline{A}}^{(2)}_{bcd} = \cyclic_{bcd}
U^{(1)p}\bas{b}(Q^{(1)}_{pcd}+\frac12 U^{(1)}_{pcd})+\frac12 U^{(1)p}
U^{(1)}_{pbcd} +\del^3_{bcd}U^{(2)}+Q^{(2)}_{bcd}
\]
where
\[
U^{(1)}_{p_1\ldots p_k}=\del^k_{p_1\ldots p_k} U^{(1)} \qquad
{U^{(1)p}}_{q_1\ldots q_k} = \del^{k+1}_{qq_1\ldots q_k}U^{(1)}\omega^{qp}
\qquad
\omega^{pq}\omega_{ql}=\delta^p_l.
\]
\end{lem}

\begin{pf}
At order 2, since $b^{(0)}=b^{(1)} = 0$, $u^{(0)}=u^{(1)}=0$,
$r^{(0)}=r^{(1)}=0$
\begin{enumerate}
\item $db^{(2)}=0$ by (\ref{eqderb}), so $b^{(2)}$ is a constant;
\item $du^{(2)}=b^{(2)}\omega$ by (\ref{eqderu});
so $b^{(2)}=0$ and $\n^0 u^{(2)}=0$;
\item the equation (\ref{eqderRicci}) at order $2$ yields that
$\del_a(r^{(2)}(\del_b,\del_c))$ is a constant;
again this implies  $u^{(2)}=0$ and $r^{(2)}(\del_b,\del_c)=a^{(2)}_{ab}$
is a constant.
\end{enumerate}

The definition of the (formal) Ricci tensor yields
$a^{(2)}_{ab} = -\del_q {A^{(2)}}^q_{ab}+{A^{(1)}}^p_{qb}{A^{(1)}}^q_{ap}$;
\newcommand{\g}[3]{{A^{(#1)}}^{#2}_{#3}}
\newcommand{\q}[3]{{Q^{(#1)}}\haut{#2}\bas{#3}}
\renewcommand{\u}[3]{{U^{(#1)}}\haut{#2}\bas{#3}}
\renewcommand{\k}[1]{k^{(2)}_{#1}}
Using lemma \ref{order1} with $\q1p{qb}={Q^{(1)}}_{qbk}\omega^{kp}$:
\[
{A^{(1)}}^p_{qb}{A^{(1)}}^q_{ap}=
\q1p{qb}\q1q{ap} +
\del_q(\q1q{ap}\u1pb)
+ \del_p(\u1qa\q1p{qb}) + \del_q(\u1pb\u1q{ap}).
\]
Hence:
\[
  a^{(2)}_{ab} =
\q1p{qb}\q1q{ap}-\del_q(\g2q{ab}-\u1pb\q1q{ap}-\u1pa\q1q{pb}-\u1pb\u1q{ap}).
\]
Since there are no exact, non-zero, \invariant\ $2n$-form on $T^{2n}$,
we have
\[
  a^{(2)}_{ab} = \q1p{qb}\q1q{ap}, \qquad
\del_q(\g2q{ab}-\u1pb\q1q{ap}-\u1pa\q1q{pb}-\u1pb\u1q{ap} )=0.
\]
The definition of the (formal) curvature tensor at order $2$ gives
$R^{(2)}_{abcd}=\del_a {\underline{A}}^{(2)}_{bcd}-
\del_b {\underline{A}}^{(2)}_{acd}
+ {A^{(1)}}^p_{bc}{\underline{A}}^{(1)}_{apd} -
{A^{(1)}}^p_{ac}{\underline{A}}^{(1)}_{bpd}.
$
Using lemma 4  we get
\begin{eqnarray*}
  R^{(2)}_{abcd} &=&
  \del_a({\underline{A}}^{(2)}_{bcd}+\u1{}{pd}\q1p{bc}-\u1pc\q1{}{bpd}
-\u1pc\u1{}{bpd}) \\
& & - \ \del_b({\underline{A}}^{(2)}_{acd}+\u1{}{pd}\q1p{ac}-\u1pc\q1{}{apd}
-\u1pc\u1{}{apd}) \\
& & +\ \q1p{bc}\q1{}{apd}-\q1p{ac}\q1{}{bpd}.
\end{eqnarray*}
The $W^{(2)}=0$ condition says that:
\begin{eqnarray*}
R^{(2)}_{abcd} &=& -\frac1{2(n+1)}\left[2\omega_{ab}a^{(2)}_{cd}
+ \omega_{ac}a^{(2)}_{bd}+ \omega_{ad}a^{(2)}_{bc}-
\omega_{bc}a^{(2)}_{ad}- \omega_{bd}a^{(2)}_{ac} \right].
\end{eqnarray*}
The fact that there does not exist a non-zero \invariant\ exact
2-form implies on one hand:
\begin{eqnarray*}
& & \del_a({\underline{A}}^{(2)}_{bcd}+\u1{}{pd}\q1p{bc}-\u1pc\q1{}{bpd}
-\u1pc\u1{}{bpd}) \\
& -& \del_b({\underline{A}}^{(2)}_{acd}+\u1{}{pd}\q1p{ac}-\u1pc\q1{}{apd}
-\u1pc\u1{}{apd}) = 0,
\end{eqnarray*}
and on the other hand:
\begin{eqnarray*}
  \q1p{bc}\q1{}{apd}-\q1p{ac}\q1{}{bpd} &=& -\frac1{2(n+1)}
\left[2\omega_{ab}a^{(2)}_{cd}
+ \omega_{ac}a^{(2)}_{bd}+ \omega_{ad}a^{(2)}_{bc}\right. \\
& & \qquad \left.-\
\omega_{bc}a^{(2)}_{ad}- \omega_{bd}a^{(2)}_{ac} \right],
\end{eqnarray*}
where $a^{(2)}_{ab}=\q1p{qb}\q1q{ap}\ $.

This last relation tells us that the \invariant\ connection
defined by $\n^0 + tQ^{(1)}$ (which is symplectic because of the complete
symmetry) has a $W$ tensor which is zero. Lifting everything to
$\R^{2n}$ and applying lemma 3 we get that the corresponding curvature
vanishes identically. Hence:
\[
  a^{(2)}_{ab}=0, \qquad \q1p{bc}\q1{}{apd}-\q1p{ac}\q1{}{bpd}=0.
\]
This in turn implies
\[
  r^{(2)} = 0, \qquad R^{(2)} = 0.
\]
The first relation tells us that there exist functions $k^{'(2)}_{cd}$
and constants $\q2{}{bcd}$ such that
\[
  {\underline{A}}^{(2)}_{bcd}-\u1pc\q1{}{bpd}-\u1{p}{d}\q1{}{bpc}
-\u1pc\u1{}{bpd} = \del_b k^{'(2)}_{cd}+\q2{}{bcd}.
\]
This can be rewritten as
\begin{equation}\label{A2}
{\underline{A}}^{(2)}_{bcd}
- \cyclic_{bcd}\u1pb(\q1{}{pcd}+\frac12\u1{}{pcd})
-\frac12\u1p{}\u1{}{pbcd}=\del_b k^{(2)}_{cd}+\q2{}{bcd}
\end{equation}
with
\[
k^{(2)}_{cd}=k^{'(2)}_{cd}-\u1p{}\q1{}{pcd}+\frac12\u1pc\u1{}{pd}
-\frac12\u1p{}\u1{}{pcd}.
\]
Indeed we have
$\u1pc\u1{}{bpd}=\frac12\u1pc\u1{}{bpd}+\frac12\del_b(\u1pc\u1{}{pd})
+\frac12\u1pd\u1{}{bpc}$ and also
$\frac12\u1pb\u1{}{cpd}=\frac12\del_b(\u1p{}\u1{}{cpd})
-\frac12\u1p{}\del_b\u1{}{cpd}$.

Now the left hand side of the equation \ref{A2} is totally symmetric
in its indices ($bcd$) so the same reasoning as in Lemma \ref{order1}
shows that $Q^{(2)}$ is totally symmetric and there exists a function
$U^{(2)}$ so that $\del_b k^{(2)}_{cd}=\del^3_{bcd}U^{(2)}$.
Substituting, we find:
\[
{\underline{A}}^{(2)}_{bcd}
= \cyclic_{bcd}\u1pb(\q1{}{pcd}+\frac12\u1{}{pcd}) +
\frac12\u1p{}\u1{}{pbcd}
+\del^3_{bcd}U^{(2)}+\q2{}{bcd}
\]
which ends the proof of the lemma.
\end{pf}

\section{A Recurrence Lemma}\label{section:recurrence}
{
\newcommand{\g}[3]{{A^{(#1)}}^{#2}_{#3}}
\newcommand{\q}[3]{{Q^{(#1)}}\haut{#2}\bas{#3}}
\renewcommand{\u}[3]{{U^{(#1)}}\haut{#2}\bas{#3}}
\renewcommand{\k}[1]{k^{(2)}_{#1}}

\begin{lem}\label{induction}
Let $\n^t$ be a formal curve of symplectic connections on $(T^{2n},\omega)$
such
that $\n^{(0)}=\n^0$, and $W^t=0$. Assume that, for all orders $l<k$,
  ${\underline{A}}^{(l)}$, and thus $r^{(l)}$, $u^{(l)}$, $b^{(l)}$ are
\invariant. Then, at order k, $r^{(k)}$, $u^{(k)}$, $b^{(k)}$ are
\invariant, and there exist a function $U^{(k)}$ on $T^{2n}$ and a
\invariant\ completely symmetric 3 tensor $Q^{(k)}$ such that
\[
{\underline{A}}^{(k)} = \del^3 U^{(k)} + Q^{(k)}.
\]
\end{lem}

\begin{pf}
Assume that, up to order
$k-1$ (included), ${\underline{A}}^{(l)}_{abc}$, $r^{(l)}_{ab}$, $u^{(l)}_{a}$,
$b^{(l)}$ are \invariant. Then, at order $k$, we have

\newcommand{\sumss}{\multiindex{\sum}{s+s'=k}{s,s'>0}}
\newcommand{\sumsss}{\multiindex{\sum}{s+s'+s''=k}{s,s',s''>0}}

\begin{eqnarray*}
&(i)& R^{(k)}_{abcd}=\del_a{\underline{A}}^{(k)}_{bcd}
-\del_b{\underline{A}}^{(k)}_{acd}
+ \sumss
\g sp{bc} {\underline{A}}^{s'}_{apd}-\g
sp{ac}{{\underline{A}}^{(s')}}_{bpd}; \\
&(ii)& r^{(k)}_{ac} = -\del_q \g kq{ac}+\sumss \g sp{qc}
\g{s'}q{ap}; \\
&(iii)& \del_c r^{(k)}_{ab}-\sumss \g sp{ca}r^{(s')}_{pb}+
{\Gamma^{(s)}}^p_{cb}r^{(s')}_{ap} =
\frac1{2n+1}(\omega_{cb}u^{(k)}_a+\omega_{ca}u^{(k)}_b); \\
&(iv)& \del_b u^{(k)}_{a} - \sumss \g sp{ba} u^{(s')}_p =
-\frac{1+2n}{2(1+n)} \sumss r^{(s)}_{bc} {r^{(s')}}^c\bas{a}
+ b^{(k)}\omega_{ba}; \\
&(v)& \del_a b^{(k)} = \frac1{1+n}\sumss \overline u^{(s)c}r^{(s')}_{ca}.
\end{eqnarray*}

Relation $(v)$ implies that $db^{(k)}$ is \invariant. Hence
$db^{(k)}=0$ and $b^{(k)}$ is a constant. Antisymmetrising $(iv)$ we
get that $ du^{(k)}-b^{(k)}\omega $ is a \invariant\ 2-form,
hence $ du^{(k)}=0$ and
\[
b^{(k)}\omega_{ba} -\frac{1+2n}{2(1+n)}
\sumss r^{(s)}_{bc}r^{(s')c}\bas a = 0.
\]
Also
\[
\del_b u^{(k)}_{a} = \sumss \g sp{ba}u^{(s')}_p.
\]
Using periodicity again and the fact that the right hand side is a
constant, we see that the $u^{(k)}_a$ are constants. Relation $(iii)$
tells us, for the same reason, that the $r^{(k)}_{ab}$ are constants.
Finally from $(i)$ and the $W^t=0$ condition, we get that
$\del_a{\underline{A}}^{(k)}_{bcd}- \del_b{\underline{A}}^{(k)}_{acd}$
is a constant hence
\begin{equation}\label{Ak}
  \del_a{\underline{A}}^{(k)}_{bcd}- \del_b{\underline{A}}^{(k)}_{acd} = 0.
\end{equation}
The reasoning of Lemma \ref{order1} applies to equation (\ref{Ak})
so  there  exist a function $U^{(k)}$ on $T^{2n}$ and a
\invariant\ completely symmetric 3 tensor $Q^{(k)}$ such that
\[
{\underline{A}}^{(k)} = \del^3 U^{(k)} + Q^{(k)}.
\]
\end{pf}

\noindent We can now proceed to the proof of the main theorem.

\begin{thm}
Let $\n^t$ be a formal curve of symplectic connections on
$(T^{2n},\omega)$ with $\n^0$ the standard connection, and $W^t=0$. Then there
exists a formal curve of symplectomorphisms $\psi_t$ such that
${\widetilde\n}^t \parnot \psi_t.\n^t$ is a formal curve of symplectic
connections which is \invariant\ and has ${\widetilde W}^t=0$, hence
is flat. In particular, $\n^t$ is flat.
\end{thm}

\begin{pf}
If $\n^t=\n^0+\sum_{k=0}^\infty t^p A^{(p)}$ is any
formal curve of symplectic connections, one defines as in \ref{actconn}
the action  of a formal curve $\psi_t$ of symplectomorphisms on $\n^t$:
\[
  (\psi_t\cdot\n^t)_X Y =
\psi_t\cdot\left(\n^t_{\psi_t^{-1}\cdot X}\psi_t^{-1}\cdot Y\right).
\]
Consider a formal one-parameter group $\psi_f(t)$ of
symplectomorphisms generated by a hamiltonian vector field $X_{f}$
($i(X_{f})\omega=df$) and consider the formal curve of symplectomorphisms
defined by
$\psi^k_f(t)=\psi_f(t^k)$.
Write
\[
  \psi^k_f(t)\cdot\n^t = \n^0+\sum_{p=0}^\infty t^p \widetilde A^{(p)}
\]
then $\widetilde A^{(p)} = A^{(p)}, \forall p<k$ and
\[
\widetilde A^{(k)}_X Y = A^{(k)}_X Y+ [X_{f},\n^0_Y Z] -
\n^0_{[X_{f},Y]} Z - \n^0_Y[X_{f},Z].
\]
Observe that $[X_{f},\n^0_Y Z] - \n^0_{[X_{f},Y]} Z - \n^0_Y[X_{f},Z] =
R^0(X_f,Y)Z+((\n^0)^2X_f)(Y,Z)$ and $\omega(((\n^0)^2X_f)(Y,Z),T) =
((\n^0)^3 f)(Y,Z,T)$.

\medskip

Assume now that the curve $\n_t=\n^0+\sum_{k=0}^\infty t^p A^{(p)}$
is a curve of symplectic connections on the torus $(T^{2n},\omega)$
and that $\n^0$ is the standard flat connection.

At order $1$, we have seen in Lemma \ref{order1} that
${\underline{A}}^{(1)} = ({\n}^0)^3U^{(1)}+ Q^{(1)}$
so choosing $f_1=-U^{(1)}$ and $\psi^{(1)}(t)=\psi_{f_1}(t)$
  as defined above we see that
\[
\psi^{(1)}(t)\cdot\n^t = \n^0+ t {\overline{Q}}^{(1)}
   +\sum_{p=2}^\infty t^p \widetilde A^{(p)}
\]
with $\omega({\overline{Q}}^{(1)}(X)Y,Z)=Q^{(1)}(X,Y,Z)$.

\medskip

Assume now that one has found a formal curve of symplectomorphisms
$ \psi^{(k-1)}(t)$ so that
\[
  \psi^{(k-1)}(t)\cdot\n^t = \n^0+ \sum_{p=1}^{k-1}t^p {\overline{Q}}^{(p)}
   +\sum_{p=k}^\infty t^p \widetilde A^{(p)}
\]
where the ${\overline{Q}}^{(p)}$ are \invariant.

At order $k$, we have seen in Lemma \ref{induction} that
${\underline{A}}^{(k)} = ({\n}^0)^3U^{(k)}+ Q^{(k)}$
where $Q^{(k)}$ is \invariant,
so choosing $f_k=-U^{(k)},~\psi^k_{f_k}(t)$ as defined above
and $\psi^{(k)}(t)= \psi_{f_k}(t^k)\circ \psi^{(k-1)}(t)$ we see that
\[
\psi^{(k)}(t)\cdot \n^t=\psi_{f_k}(t^k)\cdot\psi^{(k-1)}(t)\cdot \n^t =
  \n^0+ \sum_{p=1}^{k}t^p {\overline{Q}}^{(p)}
   +\sum_{p=k+1}^\infty t^p \widetilde A^{(p)}
\]
with $\omega({\overline{Q}}^{(k)}(X)Y,Z)=Q^{(k)}(X,Y,Z)$.
By induction this proves that one can build
a formal curve of symplectomorphisms
\[
  \psi(t)= \ldots \circ
\psi_{(f_k)}(t^k) \circ\ldots\circ \psi_{f_2}(t^2)\circ\psi_{f_1}(t)
\]
so that ${\widetilde\n}(t) \parnot \psi(t).\n(t)$ is a formal curve of
symplectic
connections which is \invariant\ and has ${\widetilde W}(t)=0$. Lifting
the connection to $\R^{2n}$ and using Lemma \ref{flatformal} shows that
${\widetilde\n}(t)$ has vanishing curvature.
Since $\n(t)=(\psi(t))^{-1}\cdot{\widetilde\n}(t)$, its curvature
is $0$ so $\n(t)$ is flat.
\end{pf}
}

The above theorem implies:
\begin{thm}
Let $\n^t$ be an analytic curve of analytic symplectic connections on
$(T^{2n},\omega)$ such that $\n^0$ is the standard flat connection
on $T^{2n}$, and such that $W^t=0$. Then
the  curvature $R^t$ of $\n^t$ vanishes.
\end{thm}

\section{Equivalence of formal curves of connections}\label{section:eq}

In this section we study the question of when two formal curves of flat
invariant connections on $T^{2n}$ are equivalent by a formal curve of
symplectomorphisms. First we consider the question on $(\R^{2n},
\Omega)$. Here it is easy to answer.

The first case to consider is the case of a single flat invariant
connection $\n^A = \n^0 +A$ on $(\R^{2n}, \Omega)$. We have seen that
such a connection is given by a linear map $A \colon \R^{2n} \to
\sp(2n,\R)$ satisfying $A(X)A(Y)=0$ and $\Omega(A(X)Y,Z)$ completely
symmetric. Define $\psi^A \colon \R^{2n} \to \R^{2n}$ by
\[
\psi^A(x) = x - \frac12 A(x)x.
\]

\begin{prop}
$\psi^A$ is a symplectomorphism of $(\R^{2n}, \Omega)$
satisfying $\psi^A \cdot \n^0 = \n^A$.
\end{prop}

\begin{pf}
It is enough to check that $\psi^A$ is a symplectomorphism on constant
vector fields. We make extensive use of the fact that $A(X)A(Y)=0$. If
$X$ is a constant vector field then
\[
\psi^A_* X_x
= \left.\frac{d}{dt} \psi^A(x+tX)\right|_{t=0}
= (X - A(x)X)_{\psi^A(x)},
\]
thus $\psi^A \cdot X = X - A(\cdot)X$. Hence
\[
\Omega(\psi^A \cdot X, \psi^A \cdot Y) (x)
= \Omega(X - A(x)X, Y - A(x)Y) = \Omega(X,Y).
\]
It is easy to see that $\psi^{-A}$ is an inverse for $\psi^A$ so that
$\psi^A$ is a symplectomorphism. Indeed, $t \mapsto \psi^{tA}$ is a
1-parameter group of symplectomorphisms with generator the symplectic
vector field $(X_A)x = -\frac12 A(x)x_x$.

Finally, for constant vector fields $X$, $Y$
\[
(\psi^A \cdot \n^0)_XY
= \psi^A \cdot (\n^0_{\psi^{-A} \cdot X}\psi^{-A} \cdot Y)
= \psi^A \cdot (  (X + A(\cdot)X) (A(\cdot)Y) ).
\]
But
\[
(X + A(\cdot)X) (A(\cdot)Y)_x
= \left.\frac{d}{dt} A(x + t(X + A(x)X))Y \right|_{t=0}
= A(X)Y
\]
so
\[
(\psi^A \cdot \n^0)_XY
=  \psi^A \cdot (A(X)Y) = A(X)Y = \n^A_XY.
\]
\end{pf}

If $\n^t = \n^0 + A^t$ is a formal curve of invariant flat connections
on $(\R^{2n}, \Omega)$ given by a curve of linear maps
$A^t \colon \R^{2n} \to \sp(2n,\R)\Bbt$
satisfying $A^t(X)A^t(Y)=0$ and $\Omega(A^t(X)Y,Z)$ completely
symmetric, we define a formal curve of vector fields $X_{A^t}$ by
\[
X_{A^t}(f) (x) = -\frac12 (A_t(x)x)_xf
\]
and set
\[
\psi_{A^t} = \exp X_{A^t}.
\]

\begin{prop}
$\psi_{A^t}$ is a formal curve of symplectomorphisms of $(\R^{2n}, \Omega)$
and $\psi_{A^t} \cdot \n^0 = \n^{A^t}$.
\end{prop}

\begin{pf}
As the exponential of a derivation, $\psi_{A^t}$ is invertible with
inverse $\exp -X_{A^t} = \psi_{-A^t}$. Moreover $\psi_{A^t}\cdot X =
\exp \ad X_{A^t} X$ and it is easy to verify that $\ad X_{A^t} X =
A^t(\cdot)X$, $(\ad X_{A^t})^2 X = 0$ so that $\psi_{A^t}\cdot X = X -
A^t(\cdot)X$ as before. Likewise $\psi_{-A^t}\cdot X = X + A^t(\cdot)X$
so that
\[
(\psi_{A^t}\cdot\n^0)_X Y
= \psi_{A^t}\cdot ( \n^0_{\psi_{-A^t}\cdot X} (Y + A^t(\cdot)Y))
= A^t(X)Y.
\]
\end{pf}

In particular the above proves
\begin{thm}
For two curves ${\widetilde{\n^t}}$ and ${\widetilde{{\n'}^t}}$
of invariant flat connections
of Ricci-type on $(\R^{2n}, \Omega)$ with
${\widetilde{\n^0}}={\widetilde{{\n'}^0}}$ the trivial connection,
there always exists a formal curve of symplectomorphisms
${\widetilde{\psi_t}}$ so that ${\widetilde{\psi_t}}
\cdot {\widetilde{\n^t}}={\widetilde{{\n'}^t}}$.
\end{thm}

Finally, we need to know what is the general form of a formal curve
of symplectomorphisms of $(\R^{2n}, \Omega)$
which fixes the trivial connection $\n^0$.

\begin{prop}
Let $\psi_t = \sigma^* \circ \exp X_t$ be a formal curve of
symplectomorphisms with $\psi_t \cdot \n^0 = \n^0$ then $\sigma(x) =Cx
+d$ and $(X_t)_x = (C_t(x) + d_t)_x$ where $C \in Sp(2n,\R),~d\in \R^{2n},
 ~C_t \in t\sp(2n,\R)\Bbt$ and $d_t \in t\R^{2n}\Bbt$.
\end{prop}

\begin{pf}
Evaluation at $t=0$ shows that $\sigma \cdot \n^0 =\n^0$ so that $\sigma
(x) = Cx+d$ where $C \in Sp(2n,\R)$ and $d\in \R^{2n}$.
 Hence $\exp X_t \cdot \n^0 =\n^0$.
$\n^0$ is the connection for which constant vector fields are parallel,
so $(\exp X_t \cdot \n^0)_X Y = 0$ for constant vector fields $X$, $Y$.
Hence $\n^0_{\exp -X_t \cdot X}\exp -X_t \cdot Y = 0$ and so $\n^0_X
\exp -X_t \cdot Y = 0$. But the only parallel vector fields for $\n^0$
are the constant fields, so $\exp -X_t \cdot Y$ is constant. The leading
term is $-t[X^{(1)},Y]$ and hence $[X^{(1)},Y]$ is constant. Since
$X^{(1)}$ is symplectic, this means $X^{(1)}_x = (C_1x+d_1)_x$ where
$C_1 \in \sp(2n,\R)$. Further $\exp tX^{(1)}$ preserves $\n^0$ and $\exp
-tX{(1)} \circ \exp X_t = \exp X'_t$ with $X'_t = O(t^2)$ so we can
recurse to conclude that $(X_t)_x = (C_t(x) + d_t)_x$ for formal curves
$C_t \in t\sp(2n,\R)\Bbt$ and $d_t \in t\R^{2n}\Bbt$.
\end{pf}

\begin{thm}
Let $\n^t$ and ${\n'}^t$ be two curves of invariant flat connections on
$T^{2n}$ with $\n^0={\n'}^0$ the trivial connection and suppose that
there is a formal curve of symplectomorphisms $\psi_t$ with $\psi_t
\cdot \n^t={\n'}^t$ then there is an element $C \in Sp(2n,\Z)$ such
that as a symplectomorphism of $T^{2n}$ we have ${\n'}^t = C \cdot \n^t$
\end{thm}

\begin{pf}
We lift the connections and $\psi_t$ to $\R^{2n}$ and denote the lifts
by a tilde. $\widetilde{\psi}_t \cdot \widetilde{\n}^t =
\widetilde{\n'}^t$. Then $\widetilde{\n}^t = \n^0 + A^t$,
$\widetilde{\n'}^t = \n^0 + B^t$ where $A^t, B^t \colon \R^{2n} \to
\sp(2n,\R)\Bbt$ are linear with the usual properties. Thus
\[
(\widetilde{\psi}_t \circ \psi_{A^t}) \cdot \n^0
=  \psi_{B^t} \cdot \n^0
\]
and hence
\[
\widetilde{\psi}_t \circ \psi_{A^t}
= \psi_{B^t} \circ \sigma^* \circ \exp X_t
\]
where $\sigma(x) = Cx+d$ and $(X_t)_x = (C_t x + d_t)_x$.

Now $\psi_{B^t} \circ \sigma^* = \sigma^* \circ {\sigma^{-1}}^*\circ
\exp X_{B^t} \circ \sigma^* =\sigma^* \circ \exp \sigma \cdot X_{B^t}$
and
\[
(\sigma \cdot X_{B^t})_x = (X_{C\cdot B^t})_x + ((C\cdot B^t)(x)d)_x
-\frac12 ((C\cdot B^t)(d)d)_x
\]
and the last two terms are in the pronilpotent semidirect product
$t\sp(2n,\R)\Bbt + t\R^{2n}\Bbt$. We can exponentiate this equation in
the form
\[
\exp \sigma \cdot X_{B^t} = \exp X_{C\cdot B^t} \exp Z_t
\]
with $Z_t \in t\sp(2n,\R)\Bbt + t\R^{2n}\Bbt$. At order zero we see that
$\sigma$ must be the lift of $\psi^{0}$ and so must preserve the
lattice: $C \in Sp(2n,\Z)$. Then $\sigma^{-1} \circ \widetilde{\psi}_t$
descends to the torus and leads off with the identity, so is of the form
$\exp L_t$ where $L_t$ is a formal series of periodic vector fields on
$\R^{2n}$. Thus we have, combining the terms in $\exp t\sp(2n,\R)\Bbt +
t\R^{2n}\Bbt$ and renaming as $Z_t$,
\[
\exp L_t = \exp X_{C\cdot B^t} \exp Z_t \exp -X_{A^t}.
\]
Equating the coefficient of t on both sides we see that
\[
L^{(1)} = X_{C\cdot B^{(1)}} + Z^{(1)} - X_{A^{(1)}}
\]
and since linear and quadratic functions are never periodic we see
that $C\cdot B^{(1)} = A^{(1)}$, and $L^{(1)}=Z^{(1)}$ is constant.
A simple recursion (moving constant terms past $\exp X_{C\cdot B^t}$)
suffices to see that $A^t = C \cdot B^t$.
\end{pf}

So we have:

\begin{thm}
The moduli space of curves of Ricci-type symplectic connections starting
with the standard flat connection on $(T^{2n},\omega)$ under the action
of formal curves of symplectomorphisms is described by the space of
formal curves $A^t \colon \R^{2n} \to \sp(2n,\R)\Bbt$ satisfying
$A^t(X)A^t(Y)=0$ and $A^t(X)Y=A^t(Y)X$, modulo the action of
$Sp(2n,\Z)$.
\end{thm}

It is worth noting that a curve of Ricci type connections on the torus
is equivalent to the constant curve at the trivial connection when
lifted to $\R^{2n}$.

{
}

\end{document}